\documentclass[12pt]{article}

\usepackage[active]{srcltx} 

\usepackage{amsfonts}
\usepackage{amsmath}
\usepackage{epsfig,psfrag}
\usepackage{graphicx}
\usepackage{epstopdf}
\usepackage{latexsym}
\usepackage[latin1]{inputenc}
\usepackage[OT1]{fontenc}
\usepackage{color,soul}

\def\R{\mathbb R}
\def\N{\mathbb N}

\def\epsilon{\varepsilon}
\def\e{\varepsilon}

\let\e=\varepsilon
\let\vp=\varphi

\let\ol=\overline
\let\ul=\underline

\let\.=\cdot
\let\0=\emptyset

\def\O{\Omega}

\def\di{\displaystyle}
\def\ckpp{c_{KPP}}

\def\pe{principal eigenvalue}
\def\pf{principal eigenfunction}
\def\MP{maximum principle}
\def\SMP{strong maximum principle}

\def\ol{\overline}

\def\qed{\hfill$\square$\\}

\def\di{\displaystyle}
\def\trait (#1) (#2) (#3){\vrule width #1pt height #2pt depth #3pt}

\newcommand{\SE}{\setcounter{equation}{0} \section}
\newcommand{\beq}{\begin{equation}}
\newcommand{\eeq}{\end{equation}}
\newcommand{\baa}{\begin{array}}
\newcommand{\eaa}{\end{array}}
\newcommand{\ba}{\begin{eqnarray}}
\newcommand{\ea}{\end{eqnarray}}

\newtheorem{thm}{Theorem}[section]
\newtheorem{theorem}{Theorem}[section]

\newtheorem{proposition}[thm]{Proposition}
\newtheorem{lem}[thm]{Lemma}

\newtheorem{remark}[thm]{Remark}

\begin{document}

\setstcolor{red}

\title{\bf The influence of a line with fast diffusion on Fisher-KPP propagation}
\author{Henri {\sc Berestycki}$^{\hbox{a }}$,  
Jean-Michel {\sc Roquejoffre}$^{\hbox{b }}$, Luca {\sc Rossi}$^{\hbox{c }}$\\
\footnotesize{$^{\hbox{a }}$ Ecole des Hautes Etudes en Sciences Sociales}\\
\footnotesize{ CAMS, 54, bd Raspail F-75270 Paris, France}\\
\footnotesize{$^{\hbox{b }}$ Institut de Math\'ematiques de Toulouse,
Universit\'e Paul Sabatier}\\
\footnotesize{118 route de Narbonne, F-31062 Toulouse Cedex 4, France}\\
\footnotesize{$^{\hbox{c }}$Dipartimento di Matematica, 
Universit\`a degli Studi di Padova}\\
\footnotesize{
Via Trieste, 63 -
35121 Padova, Italy}\\
}
\maketitle
\centerline{ \em Dedicated in friendship to Odo Diekmann}

\begin{abstract}
\noindent  We propose here a new model to describe biological invasions in the
plane when a strong diffusion takes place on a line. We establish the main
properties of the system, and also derive the asymptotic speed of spreading in
the direction of the line. For low diffusion, the line has no effect, whereas,
past a threshold, the line enhances global diffusion in the plane and the
propagation is directed by diffusion on the line. It is shown here   that the
global asymptotic speed of spreading in the plane, in the direction of the line,
grows as the square root of the diffusion on the line.
 The model is much relevant to account for the effects of fast diffusion lines
such as roads on spreading of invasive species.

\end{abstract}



\medskip


\SE{The model}

It has long been known that fast diffusion on roads can have a driving effect on the spread of epidemics. A classical example is the spread of the ``Black death" plague in the middle of  the 14th century, considered to be one of the most devastating in human history. This pandemics is known to have spread first along the silk road. After reaching the port of  Marseilles, carried by merchant boats from Crimea, it spread northwards in Europe at a fast pace along the commercial roads connecting the cities that had trade fairs. It then also spread more slowly away from the roads, inland, bringing about a dramatic invasion.  See, for instance, the account by  \cite{Sig}.

More recently, it has been observed that invasive species such as the
Processionary caterpillar of the pine tree in Europe, have been moving faster
than anticipated. One plausible explanation is that enough individuals
might have been carried on further distances than usual by vehicles
travelling on roads going through infested areas. In the same vein, the invasion of  the Aedes albopictus mosquito (also known as ``Asian tiger mosquito'')
is a concern of public health in Europe, see \cite{lemonde}. The invasion by this insect is driven by roads. Rivers may accelerate the spread of plant pathologies: see for instance \cite{JB}.

Another example of the effect of lines on propagation in open space comes from
the observation of the population of wolves in  the Western Canadian Forest. 
GPS observations reported by McKenzie et  al. \cite{McK}, 2012, suggest that
wolves move and concentrate on seismic lines. These are straight lines (with a
width of about 5m) used by oil exploration companies for testing of oil
reservoirs.
A kinetic model leading to an anisotropic diffusion equation was recently proposed by Hillen and Painter \cite {HP}; see also \cite{Hillenconf}. Performing a Chapman-Enskog analysis (with a hyperbolic scaling) they derive a drift-diffusion equation where the drift term is directed along the lines, and an anisotropic diffusion tensor. They carry out  numerical simulations that show qualitative agreement with observations: the lines clearly enhance motion.

From a mathematical point of view, this raises the question of whether the inclusion of a line with fast diffusion may affect the overall invasion speed for species that usually thrive and reproduce in regions with moderate diffusion capacities. In addition, a natural aim is to quantify this effect. 
The present paper addresses these questions and derives precise answers. 

We propose here a model where the two-dimensional environment includes a line on
which fast diffusion takes place while reproduction and usual diffusion  only
occur outside this line. Thus we consider the line $\{ (x,0)\,:\, x\in \R\}$ in
the plane $\R^2$. 
For the sake of simplicity, we will refer to the plane as ``the field" and the line as ``the road".  But clearly, the model we propose here is a very general one that is relevant not only to ecology and population dynamics but also for a wide array of situations in biology. For a single species, we consider a system that combines the density of this population in the field $v(x,y,t)$ and  the density on the line $u(x,t)$. Exchanges of populations take place between the road and the field. Namely a fraction $\nu$ of individuals from the field at the road, that is, $v(x,0,t)$, join the road while a fraction $\mu$ of the population on the road, $u(x,t)$, goes in the fields.  It is assumed that the population in the field is subject to a logistic type of growth resulting in a Fisher- KPP type of reaction term $f(v)$ in the field. We assume that no such reaction is relevant on the road. The diffusion coefficient in the field is represented by $d$ and in the road by $D$.
By reasons of symmetry, it will be enough to consider the half plane $\Omega :=
\{ (x,y) \, : \, x\in\R,\ y >0 \}$.
Translating into equations all that we have just described, we are led to the
following system:
\begin{equation}
\label{Cauchy}
\begin{cases}
\partial_t u-D \partial_{xx} u= \nu v(x,0,t)-\mu u & x\in\R,\ t>0\\
\partial_t v-d\Delta v=f(v) & (x,y)\in\O,\ t>0\\
-d\partial_y v(x,0,t)=\mu u(x,t)-  \nu v(x,0,t) & x\in\R,\ t>0.
\end{cases}
\end{equation}
The various parameters are supposed constant $D\geq0$, $d$, $\mu$, $\nu >0$ and
are given.\\

System \eqref{Cauchy} incorporates the case where the population also lives
in the lower half-plane $\R\times(-\infty,0)$ and is symmetric with respect to
the $x$ axis. Indeed, since $u$ takes the contribution by $v$
from both sides of the road and split its contribution to the upper and lower
half-plane, the problem reduces to
\eqref{Cauchy} with $\nu$ replaced by $2\nu$ in the first equation and $\mu$ by
$\frac\mu2$ in the last one. This change does not affect
the linearised system and then, as we will see in the sequel, the spreading
properties of solutions. Indeed, one come back to the original system by 
simply multiplying $v$ by $2$.

Our assumptions on the reaction term $f$ are that  
$f\in C^1([0,1])$ and satisfies 
$$f(0)=f(1)=0,\qquad\forall s\in(0,1),\quad 0<f(s)\leq f'(0)s.$$
We extend it to a negative function outside $[0,1]$.

Finally, rescaling the time by a factor $1/\nu$ in \eqref{Cauchy}, we get 
a multiplicative term $1/\nu$ in front of $\nu$, $\mu$, $D$, $d$, $f$. Thus, we
assume in the rest of the
paper that $\nu=1$; this does not entail
any loss of generality.  

The quantity 
$$
\ckpp:=2\sqrt{df'(0)}
$$
is the spreading velocity in the usual KPP equation
\cite{AW}, 
$$
\label{e4.6}
u_t-du_{xx}=f(u).
$$
As we will see, the problem has the structure of a monotone system
\cite{Hirsch}. See \cite{BES}, \cite{ES} for results of front propagation in
monotone systems, with which our problem shares many features. The originality
here lies in the 1D-2D coupling.
The main result of this paper is the
\begin{theorem}\label {t1.1} (i). Spreading.
There is an {\em asymptotic speed of spreading} $c_*=c_*(\mu,d,D)>0$ such that the following is true.
Let $(u,v)$ be a solution of \eqref{Cauchy} with a nonnegative, compactly
supported initial datum $(u_0,v_0)\not\equiv(0,0)$. Then: \\
\noindent $\bullet$ for all $c>c_*$, we have $\di\lim_{t\to+\infty}\sup_{\vert
x\vert\geq
ct}(u(x,t),v(x,y,t))=(0,0)$.\\
\noindent $\bullet$ For all $c<c_*$, we have 
$\di\lim_{t\to+\infty}\inf_{\vert x\vert\leq
ct}(u(x,t),v(x,y,t))=(1/\mu,1)$.

\medskip
\noindent (ii). The spreading velocity. If $d$ and $\mu$ are fixed, and $D$ varies in $(0,+\infty)$, the following holds true. \\
\noindent $\bullet$  If $D\leq2d$, then $c_*(\mu,d,D)=\ckpp$. \\
\noindent $\bullet$ If  $D>2d$, then $c_*(\mu,d,D)>\ckpp$ and
$\di\lim_{D\to+\infty}c_*(\mu,d,D)/\sqrt{D}$ exists and is a positive real
number.
\end{theorem}
Note that, in the statement of Theorem \ref{t1.1}, the convergence holds
pointwise in $y$. A more precise study of the domains where convergence to
$(1/\mu,1)$ or 
$(0,0)$ holds is outside the scope of this paper, and will be done elsewhere. 

In the recent years, there have been a considerable number of works on propagation in heterogeneous media, and how heterogeneities may enhance, or block,
propagation. See for instance \cite{w}, \cite{bh}. See also \cite{bhbook}
for many other references.


\SE{Conservation of total population}

If $(u,v)$ is a solution of \eqref{Cauchy} with $f\equiv0$, the
quantity $\|u(\.,t)\|_{L^1(\R)}+\|v(\.,t)\|_{L^1(\O)}$ does not depend on $t$.
To see this, suppose that $u$ and $v$ decay faster than some exponential
functions at time $t=0$. Anticipating on the next sections, we assert 
 that this property still holds for $t>0$, also for the
derivatives of $u$ and $v$, owing to parabolic estimates. We can therefore
integrate by parts the first two equations in \eqref{Cauchy} and we find:
$$\|u(\.,T)\|_{L^1(\R)}-\|u(\.,0)\|_{L^1(\R)}=\int_0^T\int_{-\infty}^{+\infty}
(v(x,0,t)-\mu u(x,t))dx\,dt,$$
\[\begin{split}
\|v(\.,T)\|_{L^1(\O)}-\|v(\.,0)\|_{L^1(\O)} &=
-d\int_0^T\int_{-\infty}^{+\infty}\partial_y v(x,0,t)dx\,dt\\
&=\int_0^T\int_{-\infty}^{+\infty}(\mu u(x,t)-v(x,0,t))dx\,dt.
\end{split}\]
Whence,
$$\|u(\.,T)\|_{L^1(\R)}+\|v(\.,T)\|_{L^1(\O)}=\|u(\.,0)\|_{
L^1(\R)}+\|v(\.,0)\|_{L^1(\O)}.$$

\medskip
\noindent{\bf Biological interpretation.} Our model is consistent with the
conservation of the total population in the case of zero natality/mortality
rate. The exchange between the line and the open plane exactly compensate each
other as is natural. The only modification in the total size of the population
comes from the ``effective birth rate" represented by the term $f(v)$.
\qed


\SE{The Cauchy problem}\label{sec:Cauchy}

In this section, we derive the existence and uniqueness result for the Cauchy
problem associated with \eqref{Cauchy}. We prove it
when $\O$ is the $N+1$-dimensional half space and, for later use, with the
addition of a drift term $q\in\R^N$, $r\in\R^{N+1}$.
Namely, denoting the generic point in $\O$ by
$(x,y)\in\R^N\times(0,+\infty)$, we study the problem
\begin{equation}\label{E}
\begin{cases}
\partial_t u-D \Delta_x u-q\.\nabla_x u=  v(x,0,t)-\mu u & x\in\R^N,\ t>0\\
\partial_t v-d\Delta v-r\.\nabla v=f(v) & (x,y)\in\O,\ t>0\\
-d\partial_y v(x,0,t)=\mu u(x,t)-  v(x,0,t) & x\in\R^N,\ t>0,
\end{cases}
\end{equation}
combined with the initial condition
\begin{equation}\label{IC}
\begin{cases}
u|_{t=0}=u_0 & \text{in }\R^N\\
v|_{t=0}=v_0 & \text{in }\O.
\end{cases}
\end{equation}
We always assume that $u_0$ and $v_0$ are nonnegative, bounded and continuous.
We look for solutions satisfying \eqref{E}-\eqref{IC} in the classical sense.

The system \eqref{E} is not very standard because the coupling appears in the
Robin boundary condition for $v$. Therefore, well-posedness has to be proved. 

\begin{proposition}\label{pro:Cauchy}
The Cauchy problem \eqref{E}-\eqref{IC}
admits a unique nonnegative, bounded solution.
\end{proposition}

The proof of the existence part is given in Appendix \ref{sec:ex}.
We derive the uniqueness part by showing, more
in general, that comparison between sub and
supersolutions\footnote{A subsolution (resp.~supersolution) is a couple
satisfying the system (in the classical sense)
with the $=$ signs replaced by $\leq$ (resp.~$\geq$) signs, which is also
continuous up to time $0$.} is
preserved 
during the evolution in Problem \eqref{E}. 
\begin{proposition}\label{comparison}
Let $(\underline u,\underline v)$ and $(\overline
u,\overline v)$ be respectively a 
subsolution bounded from above and a supersolution bounded
from below of
\eqref{E} satisfying $\underline u\leq
\overline u$ and $\underline v\leq
\overline v$ at $t=0$. Then, either $\underline u<\overline u$ and
$\underline v<\overline v$ for all $t$, or there exists $T>0$ such that 
$(\underline u,\underline v)=(\overline u,\overline v)$ for $t\leq T$.
\end{proposition}

\noindent{\sc Proof.}
Let $l$ be the Lipschitz constant of $f$. Define the 
following functions:
$$(\tilde u,\tilde v):=(\underline u,\underline v)e^{-lt},\qquad
(\check u,\check v):=(\overline u,\overline v)e^{-lt}.$$
They are sub and supersolution of the problem
\begin{equation}\label{Emonotone}
\begin{cases}
\partial_t u-D \Delta_x u-q\.\nabla_x u+(\mu+l)u=  v(x,0,t) & x\in\R^N,\
t>0 \\
\partial_t v-d\Delta v-r\.\nabla v=h(t,v) & (x,y)\in\O,\ t>0\\
v(x,0,t)-d\partial_y v(x,0,t)=\mu u(x,t) & x\in\R^N,\ t>0,\\
\end{cases}
\end{equation}
where $h(t,v):=e^{-lt}f(ve^{lt})-lv$ is
nonincreasing with respect to $v$. Let $\chi:\R\to\R$ be a nonnegative smooth
function satisfying
$$\chi(0)=0,\qquad
\chi'=0\ \text{ in }[0,1],\qquad \lim_{\rho\to+\infty}\chi(\rho)=+\infty,$$
$$((N-1)(d+D)+|q|+2|r|)|\chi'|+(2d+D)|\chi''|\leq1
\ \text{ in }\R.$$
Then, for $\e>0$, set
$$\hat u(x,t):=\check u(x,t)+\e(\chi(|x|)+t+1),\qquad
\hat v(x,y,t):=\check v(x,y,t)+\mu\e(\chi(|x|)+\chi(y)+t+1).$$
Using the fact that 
$$|\nabla_x(x\mapsto\chi(|x|))|\leq|\chi'(|x|)|,\qquad
|\Delta_x(x\mapsto\chi(|x|))|\leq |\chi''(|x|)|+
(N-1)|\chi'(|x|)|,$$
one readily checks that $(\hat u,\hat v)$ is still a supersolution of
\eqref{Emonotone}. Moreover, $(\hat u,\hat v)$ is strictly above $(\tilde
u,\tilde v)$ at $t=0$.
Assume by contradiction that this property does not hold for all $t>0$. Then,
$$T:=\sup\{\tau\geq0\ :\ \tilde u\leq\hat u\text{ in }\R^N\times[0,\tau],\
\tilde v\leq\hat v\text{ in }\ol\O\times[0,\tau]\}\in[0,+\infty).$$
It follows that $\tilde u\leq\hat u$ in
$\R^N\times[0,T]$, $\tilde v\leq\hat v$ in $\ol\O\times[0,T]$. Moreover, 
the continuity of the functions and the fact that $\hat u$ and $\hat v$ tend to
$+\infty$ as the space variable goes to infinity, uniformly in time, implies
that $T>0$ and either $\hat u-\tilde
u$ or $\hat v-\tilde v$ vanish somewhere at time $T$. If
$\min_{\R^N}(\hat u-\tilde u)(\.,T)=0$ then, since for
$x\in\R^N$, $0<t\leq T$,
\[\partial_t\hat u-D \Delta_x\hat u-q\.\nabla_x\hat u+(\mu+l)\hat u\geq  
\hat v(x,0,t)\geq\tilde v(x,0,t),\]
the parabolic \SMP\ yields $\hat u=\tilde u$ in $\R^N\times[0,T]$, which is
impossible.
Thus, $\min_{\R^N}(\hat u-\tilde u)(\.,T)>0$ and $\min_{\ol\O}(\hat v-\tilde
v)(\.,T)=0$.  
The \SMP\ implies that the latter cannot be attained inside $\O$. Thus,
necessarily,
$(\hat v-\tilde v)(\xi,0,T)=0$, for some $\xi\in\R^N$. This case is ruled out
too, because it implies
$$0\geq (\hat v-\tilde v)(\xi,0,T)-d\partial_y(\hat v-\tilde v)(\xi,0,T)
\geq\mu(\hat u-\tilde u)(\xi,T)>0.
$$
We have shown that $(\hat u,\hat v)$ is above $(\tilde
u,\tilde v)$ for all $t>0$, and then $\underline u\leq\overline u$ and
$\underline v\leq\overline v$ due to the 
arbitrariness of $\e$. Suppose now that there exists a contact point between 
$\underline u$ and $\overline u$ at some time $T$. Applying the strong \MP\
to the first equation in \eqref{E} we derive $(\underline u,\underline
v)=(\overline u,\overline v)$ for $t\leq T$. If $\underline u<\overline u$ for
all $t$, then the same arguments as before show that $\underline v<\overline v$.
\hfill$\square$\\

\medskip
\noindent{\bf Biological interpretation.}  If two pairs of population densities
(``road" and ``field") are initially
ordered, the exchanges between them always occur in such a way that the order will continue to hold separately for each population at the microscopic
level. This effect is somewhat more surprising than the mass conservation we discussed earlier.  Actually, such a monotonicity property is not a priori included in the model but rather comes out as an interesting property.  \qed

\medskip
The above comparison principle immediately extends to {\em generalised}
sub and supersolutions, given by the supremum of subsolutions and the
infimum of supersolutions respectively. In Section \ref{sec:spreading}, we will
need a
comparison principle for a more general class of subsolutions, obtained as the
supremum of two subsolutions in a given set and extended outside to the
one which is larger on the boundary. In the case of a single parabolic
equation, the comparison
principle always holds for such kind of generalised subsolutions, but this is no
longer true for our system \eqref{E}. We indeed need an extra assumption.


\begin{proposition}\label{gensub}
Let $E\subset\R^N\times(0,+\infty)$ and $F\subset\O\times
(0,+\infty)$ be two open sets and let $(u_1,v_1)$, $(u_2,v_2)$ be two
subsolutions of \eqref{E} bounded from above and satisfying
$$u_1\leq u_2\quad\text{on }(\partial E)\cap(\R^N\times(0,+\infty)),\qquad
v_1\leq v_2\text{ on }(\partial F)\cap(\O\times(0,+\infty)).$$
If the functions $\underline u$, $\underline v$ defined by
$$\underline u(x,t):=\begin{cases}
                      \max(u_1(x,t),u_2(x,t)) & \text{if }(x,t)\in\overline E\\
u_2(x,t) & \text{otherwise},
                     \end{cases}$$
$$\underline v(x,y,t):=\begin{cases}
                      \max(v_1(x,y,t),v_2(x,y,t)) & \text{if }(x,y,t)\in
\overline F\\
v_2(x,y,t) & \text{otherwise},
                     \end{cases}$$
satisfy
$$\underline u(x,t)>u_2(x,t)\ \Rightarrow\ \underline v(x,0,t)\geq v_1(x,0,t),$$
$$\underline v(x,0,t)>v_2(x,0,t)\ \Rightarrow\ \underline u(x,t)\geq u_1(x,t),$$
then, any supersolution $(\overline
u,\overline v)$ of \eqref{E} bounded from below and such that $\underline u\leq
\overline u$ and $\underline v\leq
\overline v$ at $t=0$, satisfies $\underline u\leq\overline u$ and $\underline
v\leq\overline v$ for all $t>0$.
\end{proposition}

\noindent{\sc Proof.} 
Since $\overline u\geq\underline u\geq u_2$ and $\overline v\geq\underline v\geq
v_2$ at $t=0$, Proposition \ref{comparison} yields $\overline u\geq u_2$
and $\overline v\geq v_2$ for all $t$.
Assume by contradiction that $(\overline
u,\overline v)$ is below $(\underline u,\underline v)$ somewhere.
Using the same transformations as in the proof of
Proposition \ref{comparison}, with $\e$ small enough, one reduces to the case
where $(u_1,v_1)$, $(u_2,v_2)$ are subsolutions of
\eqref{Emonotone} and $(\overline u,\overline v)$ is a supersolution of
\eqref{Emonotone} which is strictly above $(u_2,v_2)$ for all $t$, is
strictly above $(\underline u,\underline v)$ for $t$ less than some $T>0$ and
touches it somewhere at $t=T$.
We argue differently depending on the fact that $\overline u$ touches
$\underline u$ or not at $t=T$.
%
%

{\it Case 1. } $\overline u(\xi,T)=\underline u(\xi,T)$ for some
$\xi\in\R^N$.\\
Since $u_2(\xi,T)<\overline u(\xi,T)=\underline u(\xi,T)=u_1(\xi,T)$, 
it follows that $(\xi,T)\in E$ and there exists $\delta>0$ such that
$u_2<u_1=\underline u$ in 
$B_\delta(\xi)\times(T-\delta,T+\delta)\subset E$.
Let us call $Q:=B_\delta(\xi)\times(T-\delta,T)$.
For $(x,t)\in Q$,
the hypothesis on $(\underline u,\underline v)$ yields
$v_1(x,0,t)\leq \underline v(x,0,t)<\overline v(x,0,t)$. Whence,
$$\partial_t u_1-D \Delta_x u_1-q\.\nabla_x u_1+(\mu+l)u_1\leq v_1(x,0,t)<
\overline v(x,0,t)\quad\text{in }Q.$$
But $u_1\leq\overline u$ on the parabolic boundary of $Q$ and then
the strong maximum principle yields $u_1=\ol u$ in $Q$, which is impossible.

{\it Case 2. } $\underline u<\overline u$ for $t=T$ and 
$\underline v(\xi,\eta,T)=\overline v(\xi,\eta,T)$ for some 
$(\xi,\eta)\in\ol\O$.\\
Since $v_2(\xi,\eta,T)<\overline v(\xi,\eta,T)=\underline
v(\xi,\eta,T)=v_1(\xi,\eta,T)$, it follows that $(\xi,\eta,T)\in \overline F$
and there exists $\delta>0$ such that $v_2<v_1$ in
$(B_\delta(\xi,\eta)\cap\O)\times(T-\delta,T+\delta)$. By the hypothesis on
$v_1$ and $v_2$, this set cannot intersect $\partial F$ and then it is
contained in $F$, because $(\xi,\eta,T)\in \overline F$. Whence, we find that
$v_2<v_1=\underline v<\overline v$ in
$Q:=(B_\delta(\xi,\eta)\cap\O)\times(T-\delta,T)$.
If $\eta>0$, the parabolic strong maximum principle implies that $v_1=\overline
v$ in $Q$, which is impossible.
Therefore, $\eta=0$ and then the hypothesis on $\underline v$ yields
$u_1(\xi,T)\leq \underline u(\xi,T)$. Consequently,
$$0\geq (\overline v-v_1)(\xi,0,T)-d\partial_y(\overline v-v_1)(\xi,0,T)
=\mu(\overline u-u_1)(\xi,T)\geq\mu(\overline u-\underline u)(\xi,T).$$
We have reached a contradiction, because we are in the case $\underline
u<\overline u$ for $t=T$. \hfill$\square$

\begin{remark}\label{rem:EF}
It is clear that in the statement of Lemma \ref{l4.2} one can require that
$(u_1,v_1)$ is a subsolution of the three equations of \eqref{E} only in $E$,
$F$ and for $(x,0,t)\in\partial F$ respectively.
\end{remark}

%


\SE{Long time behaviour}
We derive a Liouville-type result for stationary solutions
of system \eqref{E} in the absence of drift, that is
\begin{equation}\label{stationary}
\begin{cases}
-D \Delta_x U=  V(x,0)-\mu U & x\in\R^N\\
-d\Delta V=f(V) & (x,y)\in\O\\
-d\partial_y V(x,0)=\mu U(x)-V(x,0) & x\in\R^N.
\end{cases}
\end{equation}

\begin{proposition}\label{pro:Liouville}
The unique nonnegative, bounded solutions of \eqref{stationary} are
$(U,V)\equiv(0,0)$ and $(U,V)\equiv(1/\mu,1)$.
\end{proposition}

\noindent{\sc Proof.}
Let $(U,V)$ be a nonnegative, bounded solution of \eqref{stationary}. By the
last equation of \eqref{stationary}, the statement is proved if we show that
$V\equiv0$ or $1$. Hence, by the strong \MP, it is sufficient to show that if
$V>0$ in $\O$ then $V\equiv1$. We do it in three steps.

\noindent Step 1. {\em The function $V$ satisfies}
\begin{equation}\label{infV}
\forall r>0,\quad\inf_{\R^N\times[r,+\infty)}V>0.
\end{equation}
This step is the only part of the proof in which the
absence of drift in the system is required.
Take $R>0$ large enough in such a way that the \pe\ of $-\Delta$
in $B_R\subset\R^{N+1}$, under Dirichlet boundary condition, is less than
$f'(0)/d$.
Let $\vp$ be the associated principal eigenfunction. It follows that there
exists $\e_0>0$ such that for $\e\leq\e_0$, $-d\Delta(\e\vp)<f(\e\vp)$ in $B_R$.
As a consequence, for given $\hat x\in\R^N$, $\hat y>R$, the elliptic \MP\
yields
$$\forall (x,y)\in B_R(\hat x,\hat y),\quad
V(x,y)\geq\e_0\vp(x-\hat x,y-\hat y).$$ 
This proves \eqref{infV}. 

\noindent Step 2. {\em $V\geq1$.}\\
Set $m:=\inf V$ and consider $((x_n,y_n))_n$ such that $V(x_n,y_n)\to
m$. We assume by way of contradiction that $m<1$ and we distinguish two
cases.

{\it Case 1: }$(y_n)_{n\in\N}$ tends to $0$.\\
We set
$$U_n(x):=U(x+x_n),\qquad V_n(x,y):=V(x+x_n,y).$$
By standard elliptic estimates,
we have that $(U_n)_n$ and $(V_n)_n$ converge (up to subsequences) locally
uniformly to some
functions $\check U$, $\check V$ satisfying the same system \eqref{stationary}.
Furthermore, $\check V(0,0)=m=\min_{\overline\O}\check V$. Since
$\Delta\check V=-f(\check V)/d\leq0$ in the intersection of $\O$ with a
neighbourhood of $(0,0)$, the \SMP\ and
Hopf's lemma imply that either $\check V\equiv m$ or $\check
V>m$ in $\O$ and $-\partial_y\check V(0,0)<0$. The first situation cannot occur
because $f(m)>0$ by \eqref{infV}. In the second situation,
using \eqref{stationary} we derive
$\mu\check U(0)-  m<0$. In particular, 
$$\inf U\leq\inf\check U<\frac m\mu.$$
Considering now the limit $\hat U$ of
translations of $U$ by (a subsequence of) one of its minimizing sequences, we
reduce to the case where
$$\hat U(0)=\min\hat U=\inf U,\qquad
-D\Delta_x \hat U\geq m-\mu\hat U.$$
Hence, evaluating the latter at $x=0$, we get the following
contradiction:
$$0\leq D\Delta_x \hat U(0)
\leq \mu\inf U-m<0.$$

{\it Case 2: }there exists a subsequence $(y_{n_k})_{k\in\N}$ of
$(y_n)_{n\in\N}$ with distance $r>0$ from $0$.\\
It follows from \eqref{infV} that $m>0$.
The sequence of the translated functions $V_k(x,y):=V(x+x_{n_k},y+y_{n_k})$
converges (up to subsequences) locally uniformly to a function $\check V$
satisfying
$$\Delta\check V<0\ \text{ in }B_\rho,\qquad \check V(0,0)=m=\min\check V,$$
for some $0<\rho<r$. This is impossible.
 
\noindent Step 3. {\em $V\leq1$.}\\
Assume by contradiction that $M:=\sup V>1$.
Analogously to the previous step, we can reduce to the case where
$V(0,0)=\max V=M>1$ and either $(U,V)$ satisfies \eqref{stationary}
or $-\Delta V=f(V)$ in a
neighbourhood of $(0,0)$. 
Since $f(M)<0$, the latter case is impossible and, by Hopf's lemma, we
necessarily have that $-\partial_y V(0,0)>0$. Whence 
\eqref{stationary} yields
$$\sup U\geq U(0)>\frac M\mu.$$
By usual arguments, the
limit $\hat U$ of translations of $U$ by (a subsequence of) one of its
maximizing sequence attains its maximum $\sup U$ at $0$ and satisfies 
$$D\Delta_x \hat U(0)\geq\mu\sup U-M>0.$$
This is a contradiction.
\hfill$\square$\\

\begin{theorem}\label{thm:ltb}
Let $(u,v)$ be a solution of \eqref{Cauchy} with a nonnegative, not identically
equal to zero, bounded initial datum. Then,
$$\lim_{t\to+\infty}(u(x,t),v(x,y,t))\to(1/\mu,1),$$ 
locally uniformly in $(x,y)\in\ol\O$.
\end{theorem}

\noindent{\sc Proof.}
Take $R>0$ large enough in such a way that the Dirichlet \pe\ of $-\Delta$
in $B_R\subset\R^2$ is less than $f'(0)/d$. Call $\vp$ the associated
\pf. Hence, for $\e>0$ small enough, the function $\ul V$ defined by $\ul
V(x,y):=\e\vp(x,y-R-1)$ satisfies $V\leq1$ and $-d\Delta\ul V<f(\ul V)$ in
$B_R(0,R+1)$.
Extending $\ul V$ by $0$ outside $B_R(0,R+1)$, we have that $(0,\ul V)$ is a
generalised
subsolution of \eqref{Cauchy} in the sense of Proposition \ref{gensub}. 
It is then easy to check that the solution $(\ul u,\ul
v)$ of \eqref{Cauchy}
with initial datum $(0,\ul V)$ is nondecreasing in time. 
On the other hand, the pair $(\ol U,\ol V)$ defined by
$$\ol U:=\max\left(\sup_\R
u_0,\frac1\mu\sup_{\O}v_0,\frac1\mu\right),\qquad
\ol V:=\max\left(\mu\sup_\R u_0,\sup_{\O}v_0,1\right).$$
is a supersolution of \eqref{Cauchy}. It is now sufficient to use Proposition
\ref{comparison} to infer that the solution $(\ol u,\ol v)$ with initial
datum $(\ol U,\ol V)$ is nonincreasing in time.
Therefore, owing to Proposition \ref{pro:Liouville}, as $t\to+\infty$, they both
converge to the unique positive, bounded solution of \eqref{stationary}: 
$(U,V)\equiv(\frac1\mu,1)$. This convergence is locally uniform in $(x,y)$ by
parabolic estimates.
Comparing $(u,v)$ with $(\ol u,\ol v)$, we derive
$$\limsup_{t\to+\infty}(u(x,t),v(x,y,t))\leq(1/\mu,1),$$
locally uniformly in $(x,y)\in\ol\O$.
Furthermore, by Proposition \ref{comparison} we know that $v>0$ in
$\O\times\R_+$. Hence, up to choosing a smaller $\e$ in the definition
of $\ul V$ if need be, it is not restrictive to assume that $\ul
V(x,y)<v(x,y,1)$ for $(x,y)\in\O$. Thus,
$$\liminf_{t\to+\infty}(u(x,t+1),v(x,y,t+1))\geq
\lim_{t\to+\infty}(\ul u(x,t),\ul v(x,y,t))=(1/\mu,1),$$
locally uniformly in $(x,y)\in\ol\O$.
The proof is thereby achieved.
\hfill$\square$\\

\medskip
\noindent{\bf Biological interpretation.}  Even though there is no birth rate on
the ``road", the presence of a logistic term in the ``field" is sufficient to
ensure global invasion (road and
field) by the species modelled by \eqref{Cauchy}. This is of course due to
the exchange term between the field and the road. \qed
\SE{Exponential solutions}\label{sol.exp}
The upper bound for the spreading speed will follow from the existence of
exponential solutions for the linearised problem
\begin{equation}
\label{e4.1}
\begin{cases}
\partial_t u-D \partial_{xx} u=  v(x,0,t)-\mu u & x\in\R,\ t\in\R\\
\partial_t v-d\Delta v=f'(0)v & (x,y)\in\O,\ t\in\R\\
-d\partial_y v(x,0,t)=\mu u(x,t)-v(x,0,t) & x\in\R,\ t\in\R,
\end{cases}
\end{equation}
because, by the KPP hypothesis, solutions to \eqref{e4.1}
provide entire super-solutions to \eqref{Cauchy}. 
The lower bound will follow from the
existence of solutions with complex exponent.
So,
we will be looking for solutions of the form
\begin{equation}
\label{e4.2}
(u(t,x),v(t,x,y))=(e^{\alpha(x+ct)},\gamma e^{\alpha(x+ct)-\beta y})
\end{equation}
where $\alpha$ and $\gamma$ are positive constants and $\beta$ is a real
(not necessarily positive) constant. In this section, we 
 analyse the condition on $(\alpha,\beta,\gamma)$ so that $u$ and $v$ given by
\eqref{e4.2} solve \eqref{e4.1}. 
The system on $(\alpha,\beta,\gamma)$ reads
\begin{equation}
\label{e4.3}
\left\{
\begin{array}{rcl}
-D\alpha^2+c\alpha&=&\gamma-\mu\\
-d\alpha^2+c\alpha&=&f'(0)+d\beta^2\\
d\beta \gamma&=&\mu-\gamma
\end{array}
\right.
\end{equation}
The third equation  of
\eqref{e4.3} gives $\gamma=\mu/(1+d\beta)$. Hence, the first equation, with
unknown $\alpha$ and 
parameters $D$, $d$, $\beta$, $\mu$, has the roots
$$
\alpha_D^\pm(c,\beta)=\frac1{2D}\biggl(c\pm\sqrt{c^2+\frac{4\mu
dD\beta}{1+d\beta}}\biggl).
$$
For $\alpha_D^\pm(c,\beta)$ to be real, the parameter
 $\beta$ has to be chosen in 
$(-\infty,-1/d)\cup[\beta_D(c),+\infty)$ with 
$$
\beta_D(c)=-\frac{c^2}{d(c^2+4\mu D)}.
$$  
Note that if $\beta<-1/d$ we have $\gamma<0$. This case
has therefore to be discarded. Also note that
$\di\lim_{c\to+\infty}\beta_D(c)=-1/d$. Denote 
$\Gamma_{c,D}^\pm:=\{\alpha=\alpha_D^\pm(c,\beta),\ \beta\geq\beta_D(c)\}$;
in the $(\beta,\alpha)$ plane, 
$\Gamma_{c,D}:=\Gamma_{c,D}^+\cup\Gamma_{c,D}^-$ is a smooth curve with leftmost
point $(\beta_D(c),c/2D)$. See Fig. 1. 

\begin{figure}[h]
\psfrag{bD}{\tiny$\beta_D$}
\psfrag{c/2d}{\tiny$c/2d$}\psfrag{c/2D}{\tiny$c/2D$}
\psfrag{c/D}{\tiny$c/D$}
\psfrag{a}{\scriptsize$\alpha$}\psfrag{b}{\scriptsize$\beta$}
\psfrag{Gd}{\footnotesize$\Gamma_{c,d}$}\psfrag{GD}{\footnotesize$\Gamma_{c,D}$}
\psfrag{a=}{\tiny$\alpha=(c+\sqrt{c^2+4\mu D})/2D$}
\psfrag{(a,b)}{\tiny$(\alpha_\pm,\beta_\pm)$}
\psfrag{(a-,b-)}{\tiny$(\alpha_-,\beta_-)$}
\psfrag{(a+,b+)}{\tiny$(\alpha_+,\beta_+)$}
\begin{center}
\includegraphics[width=\textwidth]{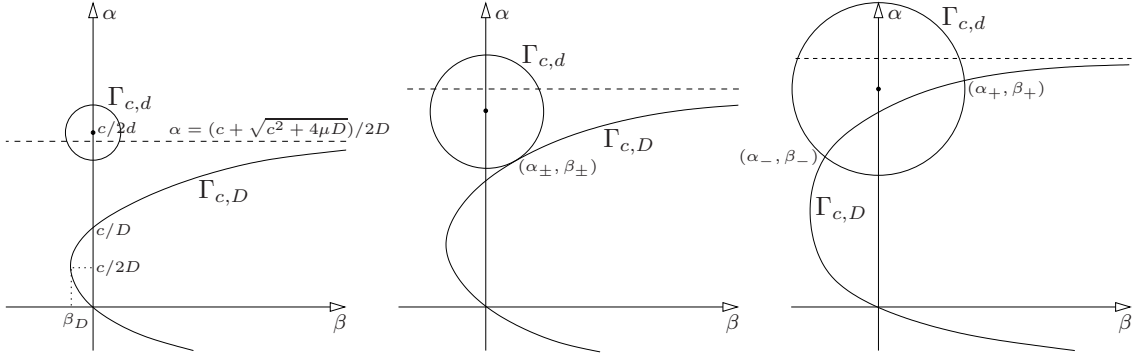}
\caption{Case  $D>2d$; $c<c_*$ (left), $c=c_*$ (middle), $c>c_*$
(right).}
\label{fig:D>2d}
\end{center}
\end{figure}

Let us now deal with the second equation in \eqref{e4.3}; in the
$(\beta,\alpha)$ plane, it is the equation of
the circle $\Gamma_{c,d}$ of centre $(0,c/2d)$ and radius
$$
\beta_{KPP}(c)=\frac{\sqrt{c^2-\ckpp^2}}{2d}.
$$
In the whole discussion below, it will be necessary to require
$$
c\geq\ckpp=2\sqrt{df'(0)}.
$$
This will be assumed without any further mention. The circle $\Gamma_{c,d}$ is
the union of  the two half circles
$\Gamma_{c,d}^{\pm}=\{\alpha=\alpha_d^\pm(c,\beta)\}$,
where
$$
\alpha_d^\pm(c,\beta)=\frac{c\pm\sqrt{c^2-\ckpp^2-4d^2\beta^2}}{2d}.
$$
 The quantity
 $$
 \alpha_d^-(c,0)=\frac{c}{2d}-\beta_{KPP}(c)
 $$
is the exponent at which a travelling wave of \eqref{e4.6} will decay
in the $x$ variable. So,
\eqref{e4.3} amounts to finding $c$ such that $\Gamma_{c,D}$
and $\Gamma_{c,d}$ intersect.
The following properties of $\alpha_D^\pm$ and $\alpha_d^\pm$ are easily
checked:
\begin{itemize}
\item The function  $\alpha_D^-$ satisfies
$\partial_\beta\alpha_D^-<0$ and $\partial_{\beta\beta}\alpha_D^->0$; moreover
$\alpha_D^-(\beta,c)<0$
if $\beta>0$, and $<0$ if $\beta\in[\beta_D(c),0)$.
\item The function 
$\alpha_D^+$ satisfies $\partial_c\alpha_D^+>0$,
$\partial_\beta\alpha_D^+>0$ wherever it is defined; moreover,
$\partial_{\beta\beta}\alpha_D^+(c,\beta)<0$.
\item The function $\alpha_d^-$ is even; moreover we have
$\partial_c\alpha_d^-<0$,  $\partial_{\beta\beta}\alpha_d^->0$, and 
$\partial_{\beta}\alpha_d^->0$ for $\beta>0$. 
\item The function $\alpha_d^+$ is even; moreover we have
$\partial_c\alpha_d^+>0$, $\partial_{\beta\beta}\alpha_d^+<0$, and 
$\partial_{\beta}\alpha_d^+<0$ for $\beta>0$. 
\end{itemize}
\noindent {\it Case 1.  $D>2d$.} For $c$ very large we have
$$\alpha_D^+(c,0)=\di\frac{c}{D}>\alpha_d^-(c,0)\sim\di\frac{f'(0)}c.
$$
So, there is $c_*>\ckpp$ such that $\Gamma_{c_*,D}$ and
$\Gamma_{c_*,d}$ are tangent. The tangency point belongs to
$\Gamma_{c_*,d}^-$; moreover, the properties $\partial_c\alpha^+_D>0$ and
$\partial_c\alpha_d^-<0$ imply that  $\Gamma_{c,D}$ and $\Gamma_{c,d}$ intersect
as soon as $c>c_*$. Due to the convexity properties of $\Gamma_{c,D}$ and
$\Gamma_{c,d}$, these two curves intersect at exactly two points denoted by
$(\beta_\pm,\alpha_\pm)$ with 
$\beta_-<\beta_+$, $\alpha_-<\alpha_+$. 
The situation is summarized in Figure \ref{fig:D>2d}.
We can be a little more explicit and tell to what part of $\Gamma_{c,D}$ the
points belong. The point here is whether $\Gamma_{c,D}^+$ and $\Gamma_{c,d}^+$
intersect, and the condition is
$$D\leq2d+\delta\quad\text{with }\delta:=
\left(\max_{t>0}\frac t{(t^2+c_{KPP}^2)(t+2)}\right)4\mu d^2.$$
Indeed $\alpha_D^+=\alpha_d^+$ admits a solution iff
$\alpha_D^+(c,\beta_{KPP}(c))\geq\alpha_d^+(c,\beta_{KPP}(c))$. The latter
reads
$$\frac{4\mu d\beta_{KPP}(c)}{c^2(1+d\beta_{KPP}(c))}\geq\frac{D-2d}{d^2},$$
whence, calling $t:=2d\beta_{KPP}(c)=\sqrt{c^2-c_{KPP}^2}$, we get
$$D-2d\leq\frac{4\mu d^2t}{(t^2+c_{KPP}^2)(t+2)},$$
which admits solutions $t>0$ iff $D\leq2d+\delta$. Actually, if
$D\leq2d+\delta$, since 
$$(D-2d)(t^2+c_{KPP}^2)(t+2)=4\mu d^2t$$
has two positive solutions (coinciding if $D=2d+\delta$), it follows that there
exist $c_*<\tilde c_1\leq\tilde c_2$ such that, for all $\tilde c_1\leq
c\leq\tilde c_2$, the equation $\alpha_D^+(c,\beta)=\alpha_d^+(c,\beta)$ has a
unique positive solution $\beta$. We have that $\delta<{\mu d}/{f'(0)}$ , that
$D\mapsto \tilde c_1(D)$ is increasing, 
$D\mapsto \tilde c_2(D)$ is decreasing and
$$\lim_{D\to2d^-}\tilde c_1(D)=0,\qquad\lim_{D\to2d^-}\tilde c_2(D)=+\infty.
$$

\noindent {\it Case 2.  $0<D=2d$.}  The centre of $\Gamma_{c,d}$ sits on
$\Gamma_{c,D}$, exactly at the point $(0,c/D)$. So, for all $c>\ckpp$,
$\Gamma_{c,D}$ and $\Gamma_{c,d}$ intersect exactly twice. In this case, we set
$c_*:=\ckpp$.

\begin{figure}[h]
\psfrag{c/2d}{\tiny$c/2d$}\psfrag{c/2D}{\tiny$c/2D$}
\psfrag{c/D}{\tiny$c/D$}
\psfrag{a}{\scriptsize$\alpha$}\psfrag{b}{\scriptsize$\beta$}
\psfrag{Gd}{\footnotesize$\Gamma_{c,d}$}\psfrag{GD}{\footnotesize$\Gamma_{c,D}$}
\begin{center}
\includegraphics[width=\textwidth]{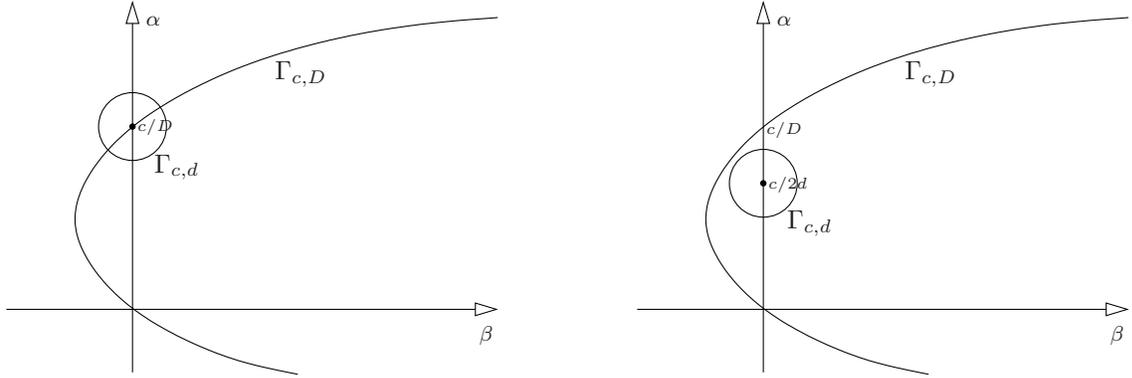}
\caption{Cases  $D=2d$ (left) and $D<2d$ (right).}
\label{fig:D<2d}
\end{center}
\end{figure}

\noindent {\it Case 3.  $0\leq D<2d$.}  Similarly to the  discussion of Case 1,
there is a largest $c'>\ckpp$ such that, for all $c<c'$, $\Gamma_{c,d}$ and
$\Gamma_{c,D}$ do not intersect. The situation
here
is more intricate than in Case 1, and also less interesting as far as spreading
is concerned: we will indeed see that here, spreading occurs at velocity
$\ckpp$. This study is therefore left for future investigation.
Setting $c_*:=\ckpp$, we see that for $c\geq c_{KPP}$, the centre
of $\Gamma_{c,d}$ lies on the right of $\Gamma_{c,D}$. This shows that the
exponential function associated with $\alpha=c/2d$, $\beta=0$ is a
supersolution of \eqref{e4.1}.\\
Cases 2 and 3 are illustrated in Figure \ref{fig:D<2d}.

\medskip
\noindent{\bf Biological interpretation.}  There is a threshold value for $D$ here, at $D=2d$. At first sight, this might appear somewhat surprising as 
one might have expected for instance $D=d$ to be the threshold (identical diffusion on the road and on the field). It turns out that this threshold value is directly related to the fact that in our model there is no reaction on the  road.
As a matter of  fact we may  consider a related model where one allows reproduction on the road with the same rate as in the field. In this case, the new  system reads:
\begin{equation} 
\label{entire1}
\begin{cases}
\partial_t u-D \partial_{xx} u= \nu v(x,0,t)-\mu u +f(u)& x\in\R,\ t\in\R\\
\partial_t v-d\Delta v=f(v) & (x,y)\in\O,\ t\in\R\\
-d\partial_y v(x,0,t)=\mu u(x,t)-  \nu v(x,0,t) & x\in\R,\ t\in\R.
\end{cases}
\end{equation}
Adapting the arguments of the next section one can show that then,
the threshold occurs at $D=d$. A detailed study of the system with a reaction 
term on the road too will be done elsewhere.
 \qed

\SE{Asymptotic spreading}\label{sec:spreading}
Because of the comparison principle (Proposition \ref{comparison}), the standard
way to prove that solutions spread at least at speed $c_*$ is to find a
compactly supported generalised stationary subsolution
in the moving framework at 
velocity $c<c_*$.
Here, generalised subsolution is in the sense of Proposition \ref{gensub}.
We consider the linearised system
penalised by $\delta>0$ in the moving framework:
\begin{equation}
\label{e4.7}
\begin{cases}
\partial_t u-D\partial_{xx} u+c\partial_x u=  v(x,0,t)-\mu u & x\in\R,\ t\in\R\\
\partial_t v-d\Delta v+c\partial_x v=(f'(0)-\delta)v & (x,y)\in\O,\ t\in\R\\
-d\partial_y v(x,0,t)=\mu u(x)-v(x,0,t) & x\in\R,\ t\in\R.
\end{cases}
\end{equation}

For $D>2d$, the main lemma is the following.
\begin{lem}
\label{l4.1}
Assume $D>2d$ and let $c_*$ be as in Section \ref{sol.exp}. For $c<c_*$ close
enough to $c_*$, there exists $\delta>0$ such that \eqref{e4.7}
admits a nonnegative, compactly supported, generalised stationary subsolution
$(\underline u,\underline v)\not\equiv(0,0)$.
\end{lem}
\noindent{\sc Proof.} In order to keep the notation as light as possible, we
will
carry out the discussion with $f'(0)$ instead of $f'(0)-\delta$ in \eqref{e4.7},
because all the results will perturb for small $\delta>0$.  

The first step is to devise a stationary solution of \eqref{e4.7}, not in $\O$
anymore, but
in the horizontal strip $\O^L:=\R\times(0,L)$ with $L>0$. So, we are solving
\begin{equation}
\label{e4.71}
\begin{cases}
-D U''+c U'=V(x,0)-\mu U & x\in\R\\
-d\Delta V+c\partial_x V=f'(0)V & (x,y)\in\O^L\\
-d\partial_y V(x,0)=\mu U(x)-V(x,0) & x\in\R\\
V(x,L)=0 & x\in\R.
\end{cases}
\end{equation}
In a similar fashion as in the preceding section, we look for exponential
solutions of the form $(1,\gamma(y))e^{\alpha x}$. 
An easy computation shows that such a solution exists if and only if the system
with unknowns $\alpha$ and $\beta$ (look for $\gamma$ under the form
$\gamma_1e^{-\beta
y}+\gamma_2e^{\beta y}$)
\begin{equation}
\label{e4.34}
\left\{
\begin{array}{rcl}
-D\alpha^2+c\alpha+\di\frac{(1+e^{-2\beta L})d\beta\mu}{1-e^{-2\beta
L}+(1+e^{-2\beta L})d\beta}&=&0\\
-d(\alpha^2+\beta^2)+c\alpha&=&f'(0)\\
d\beta (\gamma_1-\gamma_2)&=&\mu-(\gamma_1+\gamma_2)\\
\gamma_1e^{-\beta
L}+\gamma_2e^{\beta L}&=&0
\end{array}
\right.
\end{equation}
has a solution. The first equation of \eqref{e4.34} defines, in the upper part
of the plane $(\beta,\alpha)$,  a curve $\Gamma_{c,D}^L$, which is symmetric
with respect to the 
vertical axis. Moreover, the set
$\Gamma_{c,D}^L\cap\{\alpha>0,\beta>0\}$ tends to
$\Gamma_{c,D}\cap\{\alpha>0,\beta>0\}$ as $L\to\infty$. In fact, it is the graph
$\alpha=\alpha_{D}^{+,L}(c,\beta)$ with
$$
\alpha_{D}^{+,L}(c,\beta)=\frac1{2D}\biggl(c\pm\sqrt{c^2+\frac{4(1+e^{-2\beta
L})\mu dD\beta}{1-e^{-2\beta L}+(1+e^{-2\beta L})d\beta}}\biggl).
$$
The curve $\Gamma_{c,D}^L$ is strictly above $\Gamma_{c,D}$; moreover, as
$L\to\infty$, $\alpha_{D}^{+,L}$ converges to $\alpha_{D}^+$, together with
all its derivatives,  in every compact subset in $(c,\beta)$ avoiding
$\{\beta\leq0\}$. Thus, due to the implicit function theorem, for large
$L$, the picture is analogous to the case 1 of the previous section: there
exists a
unique $c_*^L\in(c_{KPP},c_*)$ such that $\Gamma_{c,D}^L\cap\{\beta>0\}$
intersects 
$\Gamma_{c,d}$ twice if $c>c_*^L$ once if $c=c_*^L$ and never if
$c<c_*^L$. Moreover we have $\di\lim_{L\to\infty}c_*^L=c_*$. 

To obtain the compactly supported subsolution, 
fix $L>0$ suitably large and set
$$
h^L(c,\beta):=\alpha_D^{+,L}(c,\beta)-\alpha_d^-(c,\beta).
$$
Call $(\beta_*^L,\alpha_*^L)$ the tangent point between
$\Gamma_{c_*^L,D}^L\cap\{\beta>0\}$ and $\Gamma_{c_*^L,d}$.
We have
$$
\partial_\beta h^L(c_*^L,\beta_*^L)=0,\qquad
-2a:=\partial_{\beta\beta}h^L(c_*^L,\beta_*^L)<0.
$$
Also, set $e:=\partial_ch^L(c_*^L,\beta_*^L)>0$,
$b:=\partial_{c\beta}h^L(c_*^L,\beta_*^L)$, and, because we are working in a
vicinity
of $(c_*^L,\beta_*^L)$:
$$
\xi:=c_*^L-c,\qquad  \tau:=\beta-\beta_*^L.
$$
The equation $h^L(c,\beta)=0$ becomes, for $(c,\beta)$ in a neighbourhood of
$(c_*^L,\beta_*^L)$:
\begin{equation}
\label{e4.8}
a\tau^2+b\xi\tau+e\xi=\varphi(\tau,\xi)
\end{equation}
where $\varphi$ is analytic in $\tau$ in a neighbourhood of 0, vanishing at
$(0,0)$ like $\vert\tau\vert^3+\xi^2$. For small $\xi>0$, the trinomial
$a\tau^2+d\xi\tau+e\xi$ has two roots $\tau_\pm=\pm
i\sqrt{(e/a)\xi}+O(\xi)$; by (an adaptation of) Rouch\'e's theorem
(see Appendix \ref{sec:Rouche}),
equation \eqref{e4.8}
has two roots, conjugate up to the order $\xi$, still called $\tau_\pm$.
They also satisfy 
$\tau_\pm=\pm i\sqrt{(e/a)\xi}+O(\xi)$. Reverting to 
the full notations, we see that, for $c$ strictly less than and suitably close
to $c_*^L$,
equation \eqref{e4.34} has a solution $(\beta,\alpha)$ with the following
properties: $\beta$ has real part $\xi$-close to $\beta_*^L$, hence positive,
and it also has nonzero imaginary part of order
$\xi^{1/2}$. It follows from the second equation of \eqref{e4.34} that $\alpha$
has nonzero imaginary part too. Therefore, one can write
$$
(\beta,\alpha)=(\beta_1+i\beta_2,\alpha_1+i\alpha_2),
$$
with $\beta_1>0$ and $\alpha_2,\beta_2\neq0$. We thus obtain a solution
$(U(x),V(x,y))=(1,\gamma(y))e^{(\alpha_1+i\alpha_2)x}$ to \eqref{e4.71}, with
$$
\gamma(y)=\gamma_1(e^{-\beta y}-e^{\beta(-2L+y)}),\  \ 
\gamma_1=\frac{\mu}{1-e^{-2\beta L}+d\beta(1+2e^{-2\beta L})}.
$$
We have $\mathrm{arg}(\gamma_1)=O(\sqrt{\xi})$. Furthermore:

\noindent$\bullet$ $\mathrm{Re}(U)>0$  on
$(-\pi/2\alpha_2,\pi/2\alpha_2)$ and vanishes at the ends;

\noindent $\bullet$ 
$\mathrm{Re}(V)>0$ if and only if
$$
\cos(\alpha_2x+\mathrm{arg}(\gamma_1)-\beta_2y)>e^{-2\beta_1(L-y)
}\cos(\alpha_2x+\mathrm{arg}(\gamma_1)-\beta_2(2L-y)).
$$
That is,
$$
\cos(\alpha_2 x+\mathrm{arg}(\gamma_1)-\beta_2y)>\frac{\sin(2\beta_2(L-y))}
{e^{2\beta_1(L-y)}-\cos(2\beta_2(L-y))}
\sin(\alpha_2x+\mathrm{arg}(\gamma_1)-\beta_2y).$$
Thus, the set where $\mathrm{Re}(V)>0$ is periodic in the direction $e_1$, with
period
$\frac{2\pi}{\alpha_2}$. Its connected
components intersecting the strip $\R\times(0,L)$ are bounded.
Since $\beta_2=O(\sqrt{\xi})$, by decreasing $\xi$ we can make one of
them, denoted by $F$, satisfy the property that $\{\alpha_2
x+\mathrm{arg}(\gamma_1)\ :\ (x,0)\in\overline F\}$ is arbitrarily close to the
set $[-\pi/2,\pi/2]$. Since
$\mathrm{arg}(\gamma_1)=O(\sqrt{\xi})$, we have that $\{\alpha_2
x\ :\ (x,0)\in\overline F\}$ is
close to $[-\pi/2,\pi/2]$ for $\xi$ small.
We define the following functions: 
$$\underline u(x):=\begin{cases}
                      \max(\mathrm{Re}(U(x)),0) & \text{if
}|x|\leq\frac\pi{2\alpha_2}\\
0 & \text{otherwise},
                     \end{cases}$$
$$\underline v(x,y):=\begin{cases}
                      \max(\mathrm{Re}(V(x,y)),0) & \text{if }
(x,y)\in\overline F\\
0 & \text{otherwise}.
                     \end{cases}$$
The choice of $F$ implies that the couple $(\underline u,\underline v)$ is a
generalised subsolution to \eqref{e4.7} in the
sense of Proposition \ref{gensub}.
\hfill$\square$\\

Turn to the case $0\leq D\leq 2d$. The main lemma is
here the following:
\begin{lem}
\label{l4.2}
Assume $0\leq D\leq2d$. For $-\ckpp<c<\ckpp$, the conclusion of Lemma \ref{l4.1}
holds.
\end{lem}
\noindent{\sc Proof.} 
Let $|c|<\ckpp$. For $0<\delta<f'(0)-c^2/4d$, the equation
\begin{equation}
\label{e4.100}
-d\Delta V+c\partial_x V=(f'(0)-\delta/2)V,\  \  \    (x,y)\in\R^2
\end{equation}
has a compactly supported subsolution. The construction is classical: setting
$\omega:=d^{-1/2}\sqrt{f'(0)-\delta-c^2/4d}$, we see that the Dirichlet
problem
$$
-d\phi''+c\phi'=(f'(0)-\delta)\phi \quad  \hbox{in
$(-\pi/2\omega,\pi/2\omega)$},\qquad \phi(\pm \pi/2\omega)=0
$$
admits the positive solution $\phi(x)=e^{(c/2d)x}\cos(\omega x)$. Then, let
$\psi_R(y)$ be the first eigenfunction of $-\partial_{yy}$ in $(-R,R)$; choosing
$R$ such that the associated eigenvalue is $\delta/2d$, it is readily seen
that $V(x,y):=\phi(x)\psi_R(y-R-1)$ is a solution to
\eqref{e4.100} in the rectangle $(-\pi/2\omega,\pi/2\omega)\times
(1,2R+1)$, vanishing on the boundary. Extending it by $0$ outside,
we find that $(0,V)$ is a generalized subsolution to \eqref{e4.7} (cf.~Remark
\ref{rem:EF}). \hfill$\Box$
\\

\noindent{\sc Proof of Theorem \ref{t1.1}: the 'spreading' part}. 
Let $c_*$ be as in Section \ref{sol.exp} and let $(e^{\alpha(x+c_*t)},\gamma
e^{\alpha(x+c_*t)-\beta y})$ be the exponential
entire supersolution of \eqref{e4.1} - and then of \eqref{Cauchy} -
constructed there. By symmetry, $(e^{\alpha(-x+c_*t)},\gamma
e^{\alpha(-x+c_*t)-\beta y})$ is still a supersolution. 
The comparison principle - Proposition \ref{comparison} - then yields
$$\forall c>c_*,\quad\di\lim_{t\to+\infty}\sup_{\vert x\vert\geq
ct}(u(x,t),v(x,y,t))=(0,0).$$

Fix now $0<c<c_*$ (close enough to $c_*$ if $D>2d$) and consider the
pair $(\underline u,\underline v)$ given by Lemmas \ref{l4.1}, \ref{l4.2}.
There exists $\gamma_0>0$ such that, for $0<\gamma\leq\gamma_0$,
$\gamma(\underline u,\underline v)$ is a subsolution of the problem in
the moving framework:
\begin{equation}
\label{moving}
\begin{cases}
\partial_t u-D \partial_{xx} u+c\partial_x u= v(x,0,t)-\mu u & x\in\R,\
t>0\\
\partial_t v-d\Delta v+c\partial_x v=f(v) & (x,y)\in\O,\ t>0\\
-d\partial_y v(x,0,t)=\mu u(x,t)-v(x,0,t) & x\in\R,\ t>0.
\end{cases}
\end{equation}
For $0<\gamma\leq\gamma_0$, let $(u_\gamma,v_\gamma)$ be the solution of
\eqref{moving} with initial datum
$\gamma(\underline u,\underline v)$. Using the comparison principle for
generalised
subsolutions - Proposition \ref{gensub} - we see that $(u_\gamma,v_\gamma)$ is
nondecreasing in $t$. Moreover, since $\gamma(\underline
u,\underline
v)$ is not a solution of \eqref{moving}, it follows from the strong comparison
principle given by Proposition \ref{comparison} that $(u_\gamma,v_\gamma)$
is strictly increasing in $t$. Thus, as $t\to+\infty$, $(u_\gamma,v_\gamma)$
converges locally uniformly to a stationary solution $(U_\gamma,V_\gamma)$ of
\eqref{moving} which is strictly larger than $\gamma(\underline u,\underline
v)$.
There exists then $k>0$ such that $(U_\gamma,V_\gamma)$ is above the translated
by any $h\in(-k,k)$ in the $x$-direction of $\gamma(\underline u,\underline v)$.
By
comparison with the translated by $h$ of $(u_\gamma,v_\gamma)$, we infer that  
$(U_\gamma,V_\gamma)$ is above the translated by $h$ of itself, that is, it does
not depend on
$x$.
Therefore, $U_\gamma$ is constant and $V_\gamma=V_\gamma(y)$ satisfies
$$V_\gamma(0)=\mu U_\gamma,\qquad-dV_\gamma''=f(V_\gamma)\ \text{ for
}y>0,\qquad V_\gamma'(0)=0.$$
This easily implies $V_\gamma\equiv1$ and $U_\gamma\equiv
1/\mu$.
We will now conclude by comparing the solution $(u,v)$ with a combination of
the functions
$$(u_\gamma(x+ct,t),v_\gamma(x+ct,y,t)),\qquad
(u_\gamma(-x+ct,t),v_\gamma(-x+ct,y,t)),$$
which are also solutions of \eqref{Cauchy}. Notice that, up to waiting until
time $1$, it is not restrictive to assume that, at initial time, $(u,v)$ is
positive and then it is above $\gamma(\underline u,\underline v)$ for
some $\gamma\leq\gamma_0$. Take $T>2$
and $0\leq\xi\leq c(T-2)$. Let $\tau\in[1,T/2]$ be such that $\xi=c(T-2\tau)$.
By comparison, we get
$$\forall (x,y)\in\O,\ 0<t<\tau,\quad
(u(x,t),v(x,y,t))\geq (u_\gamma(x+ct,t),v_\gamma(x+ct,y,t)),$$
whence, in particular,
$$(u(x,\tau),v(x,y,\tau))\geq
(u_\gamma(x+c\tau,\tau),v_\gamma(x+c\tau,y,\tau))\geq
(u_\gamma(x+c\tau,1),v_\gamma(x+c\tau,y,1)).
$$
Taking $0<\gamma'\leq\gamma_0$ small enough (independently on $T$) in such a way
that $\gamma' 
(\underline u(-x),\underline v(-x,y))\leq(u_\gamma(x,1),v_\gamma(x,y,1))$, we
find that, at $t=\tau$, $(u,v)$ is larger than the pair
$$(u_{\gamma'}(-x+c(t-2\tau),t-\tau),v_{\gamma'}(-x+c(t-2\tau),y,t-\tau)).$$
Applying once again the comparison principle we get, for given $y>0$,
$$(u(\xi,T),v(\xi,y,T))\geq(u_{\gamma'}(0,T-\tau),v_{\gamma'}
(0,y,T-\tau))\geq(u_{\gamma'}(0,T/2),v_{\gamma'}
(0,y,T/2)).$$
From this, we eventually infer 
$$\lim_{T\to+\infty}\inf_{0\leq\xi\leq
c(T-2)}(u(\xi,T),v(\xi,y,T))\geq(1/\mu,1).$$
The reverse inequality follows from Theorem
\ref{thm:ltb}. This concludes the proof, as the
negative values of $\xi$ can be handled by reflection with respect to
the $y$ axis. 
\hfill$\square$\\

\medskip
\noindent{\bf Biological interpretation.}  The result we have derived in this section is a remarkable effect of
 this model.
Essentially, it says that propagation on the road can be accelerated by a large diffusion, provided that 
there is an exchange between the road and an outside region where reproduction and some diffusion occurs. 
 For lower values of the diffusion on the road, namely, $D< 2d$, the propagation on the road is 
purely driven by that in the field. The invasion velocity then is  the standard KPP invasion speed. However, when the critical value $D=2d$ is crossed, a transition to a new regime takes place. Namely, the {\em combination} of reaction in the field and fast diffusion on the line yields an invasion speed strictly larger than the KPP invasion speed. One may view this regime as driven by the diffusion on the line. Note that it is rather remarkable that even though no reproduction takes place there (leading to dispersion if it were an isolated line), the large diffusion drives the overall phenomenon.

The precise value of the threshold diffusion here is $D=2d$. In the previous
section, we saw that when the reaction (effective reproduction) also takes place
on the road, one is led to system~\eqref{entire1}. For that system, the
threshold value is $D=d$. 
It is interesting to note that while in both cases a large diffusion on  the  
road  always enhances the spreading speed, when there is pure diffusion and no
reproduction  on the road, this effect starts to be felt at higher values of the
diffusion when compared to the model in which  the reaction holds everywhere. 
\qed

\SE{The large diffusion limit of $c_*$}

In this section we investigate the behaviour of the critical speed $c_*$,
introduced in Section \ref{sol.exp}, as $D$ goes to $\infty$. Hence, $c_*$ is
now understood as a function of the variable $D$.
We will first show that
\begin{equation}\label{e:c*}
\sqrt{4\mu^2+(f'(0))^2}-2\mu
\leq\liminf_{D\to\infty}\frac{c_*^2}D\leq
\limsup_{D\to\infty}\frac{c_*^2}D\leq f'(0).
\end{equation}
The proof of the last statement of Theorem \ref{t1.1} will follow easily.

For $D>2d$, we know that $c_*>c_{KPP}=2\sqrt{df'(0)}$. Moreover,
it follows from geometrical
considerations (see figure 1) that
$$\frac{c_*}D<\frac{c_*}{2d}-\beta_{KPP}(c_*)<\frac1{2D}(c_*+
\sqrt{c_*^2+4\mu D}),$$
where, we recall, $\beta_{KPP}(c)=\frac{\sqrt{c^2-\ckpp^2}}{2d}$ is the radius
of the circle $\Gamma_{c,d}$. Whence,
$$\frac1D<\frac1{2d}\left(1-\sqrt{1-\frac{c_{KPP}^2}{c_*^2}}
\right)<\frac1 { 2D } \left(1+
\sqrt{1+\frac{4\mu D}{c_*^2}}\right).$$
From the second inequality above we infer that, as $D\to\infty$, $c_*$ tends to
$\infty$ and then 
$$\frac1D<\frac1{2d}c_*^{-2}\left(\frac{c_{KPP}^2}2+o(1)\right)<\frac1 { 2D }
\left(1+
\sqrt{1+\frac{4\mu D}{c_*^2}}\right).$$
Calling $\tau:=\frac{c_*}{\sqrt D }$, the first inequality above yields
$\tau^2<f'(0)+o(1)$,
and then, by the second one,
$$2f'(0)+o(1)<\tau^2+\tau\sqrt{
\tau^2+4\mu}<f'(0)+\tau\sqrt{
\tau^2+4\mu}+o(1).$$
Property \eqref{e:c*} then follows.
\\

\noindent{\sc Proof of Theorem \ref{t1.1}: the limit of $c_*$.}  Revert to the
algebraic system \eqref{e4.3}. Estimate \eqref{e:c*}  enables us to rescale the
unknowns $c$ and $\alpha$ as
$$
c=\sqrt D\tilde c,\qquad\alpha=\tilde\alpha/\sqrt D,
$$
with $\tilde c$ bounded from below away from 0 and from above independently on
$D$. It then follows from the first equation of \eqref{e4.3} that 
$\tilde\alpha$ is bounded from above too.
Drop the tildes. Because $c$ and $\alpha$ are now to be found bounded, we
drop the term $-d\alpha^2/D$ produced by the
normalisation in \eqref{e4.3}, and the limiting system to solve is
\begin{equation}
\label{e7.1}
\left\{
\begin{array}{rcl}
-\alpha^2+c\alpha&=&-\di\frac{d\beta\mu}{1+d\beta}\\
\alpha&=&\di\frac1c(f'(0)+d\beta^2).
\end{array}
\right.
\end{equation}
The last equation of \eqref{e7.1} is that of a parabola in the $(\beta,\alpha)$
plane, whose bottom is the point $(0,f'(0)/c)$.  Call it $\Gamma_{c,d}^\infty$ -
just to be consistent with the notations of Section 5. The first equation just
represents the set  
$\Gamma_{c,D}$ for $D=1$.  When $c$ is close to 0, the parabola
$\Gamma_{c,d}^\infty$ is above this set - which lies below the
line $\{\alpha=2\sqrt\mu\}$. When $c$ is close to infinity, because of their
respective convexity properties, the curves intersect at two points (see Figure
\ref{fig:Dinfty} below). Thus (same
argument as in Section 5) there is $c_*>0$ such that \eqref{e7.1} has exactly
one solution. This is the sought for limit. 
\begin{figure}[h]
\psfrag{c/2d}{\tiny$c/2d$}\psfrag{c/2D}{\tiny$c/2D$}
\psfrag{c/D}{\tiny$c/D$}
\psfrag{f'}{\tiny$f'(0)/c$}
\psfrag{a}{\scriptsize$\alpha$}\psfrag{b}{
\scriptsize$\beta$}
\psfrag{Gd}{\footnotesize$\Gamma_{c,d}^\infty$}\psfrag{GD}{\footnotesize$\Gamma_
{c,D}$}
\begin{center}
\includegraphics[width=\textwidth]{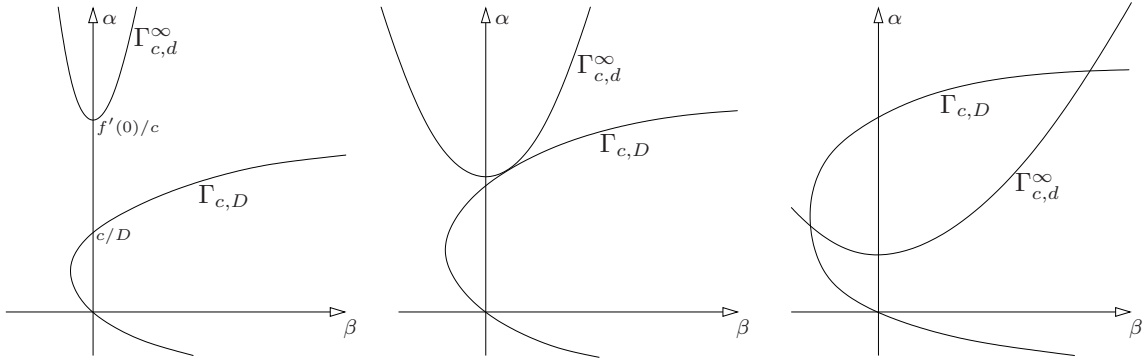}
\caption{$c$ close to $0$ (left), $c=c_*$ (middle), $c$ close to $\infty$
(right).}
\label{fig:Dinfty}
\end{center}
\end{figure}
\hfill$\square$

\medskip
\noindent{\bf Biological interpretation.}  When the diffusion on the road is very large, the spreading along the road also becomes very large. The invasion speed can be made as large as one wishes provided the diffusion on the single line is taken sufficiently large. This model therefore accounts for the enhancement of invasion speeds in the presence of lines with large diffusion coefficients even though no growth of population takes place there. 
The precise quantification we derive here states that the invasion speed grows like the square root of the pure diffusion on the road. It  seems interesting to go further into this analysis and determine how the prefactor of $\sqrt{D}$ depends on the various parameters. 
\qed


\section{Conclusion}
We have introduced a model that describes the effect of a line
on which there is fast diffusion on the overall  propagation of a species that diffuses with another constant and reproduces outside this line in a two dimensional framework.

We have found that this model conserves the population in absence of reproduction and mortality and preserves order. Then, we have shown that owing to the exchanges taking place between the line and the plane, there is an asymptotic speed of spreading which is the invasion velocity along the line. 

We have computed the global asymptotic speed of spreading along the line. This
is achieved with exponential solutions of the linearised system and compactly
supported sub-solutions. The asymptotic speed is derived from an algebraic
system. 

When $D$, the diffusion on the road, is less than or equal to $2d$, where $d$ is the diffusion in the field, there is no effect at all of the road: the propagation takes place at the classical  KPP invasion speed. In contradistinction with this case, when $D$ is larger than $2d$, there is an enhancement effect of the diffusion on the road leading to a speed higher than KPP. Lastly, this invasion speed is shown to behave like $\sqrt{D}$ for large values of $D$.



\begin{appendix}
\numberwithin{equation}{section}

\SE*{Appendix}

\SE{Existence result for the Cauchy problem}\label{sec:ex}

\noindent{\sc Proof of the existence part of Proposition \ref{pro:Cauchy}.}
We prove the result for an initial datum $(u_0,v_0)$ which is
locally H\"older continuous, together with its
derivatives up to order 2, and satisfies the compatibility condition
$$
-d\partial_y v_0(x,y,0)=\mu u_0(x,y)-  v_0(x,y,0).
$$
The regularity of the initial datum is therefore inherited by the solution of
the Cauchy problem for all time $t\geq0$. The case of a merely continuous
initial datum can then be handled by a standard regularization technique (see,
e.g., \cite{Lady}).

We will obtain a solution to \eqref{E}-\eqref{IC} as the limit of a
subsequence of
solutions $((u_n,v_n))_n$ of the following problems:
\begin{equation}\label{un}
\begin{cases}
\partial_t u_n-D \Delta_x u_n-q\.\nabla_x u_n+\mu u_n=  v_{n-1}(x,0,t)
& x\in\R^N,\ t>0\\
u_n|_{t=0}=u_0 & \text{in }\R^N,\\
\end{cases}\end{equation}
\begin{equation}\label{vn}
\begin{cases}
\partial_t v_n-d\Delta v_n-r\.\nabla v_n=f(v_n) 
& (x,y)\in\O,\ t>0\\
  v_n(x,0,t)-d\partial_y v_n(x,0,t)=\mu u_n(x,t) & x\in\R^N,\
t>0\\
v_n|_{t=0}=v_0 &  \text{in }\O,
\end{cases}\end{equation}
starting from $v_0$.

\noindent Step 1. {\em Solvability of \eqref{un}, \eqref{vn}.}\\
We say that a function $w(z,t)$ has admissible growth in $z$ if it
satisfies $|w(z,t)|\leq\beta e^{\sigma|z|^2}$, for some $\sigma,\beta>0$. It is
well known that the linear Cauchy
problem is uniquely solvable in the class of functions with admissible growth
in the space variable.
If $v_{n-1}$ is a continuous function with admissible growth, then problem
\eqref{un} admits a unique classical solution $u_n$ with admissible growth.
In order
to solve \eqref{vn}, notice that it can be reduced to a homogeneous system by
replacing $v_n$ with $v_n-v_0-\mu (u_n-u_0)$. It then follows from the
standard parabolic theory that it admits a unique classical solution with
admissible growth. Let $((u_n,v_n))_n$ denote the family of solutions
constructed in this way, starting from $v_0$.

\noindent Step 2. {\em $L^\infty$ estimates.}\\
We show, with a recursive argument, that 
$$\forall n\in\N,\quad 0\leq u_n\leq \frac1\mu H,
\quad0\leq v_n\leq H,\quad
\text{with }H:=\max\left(\mu\|u_0\|_\infty,\|v_0\|_\infty,1\right).$$
The property trivially holds for $n=0$. Assume that it holds for some
value $n-1$.
Since $0$ and $\di\frac1\mu H$ are respectively a sub and a supersolution of
\eqref{un}, the comparison principle yields $0\leq u_n\leq\frac1\mu H$.
It then follows that $H$ is a supersolution of \eqref{vn}, whence $0\leq
v_n\leq H$.

\noindent Step 3. {\em $W^{2,1}_p$ estimates.}\\
By step 1 we know that $0\leq v_{n-1}\leq H$.
Thus, applying the local boundary estimates to \eqref{un} we infer that, for any
given $\rho,T>0$ and $1<p<\infty$,
$$\|u_n\|_{W^{2,1}_p(B_{\rho+1}\times(0,T))}\leq C H,$$
where $B_\rho$ denotes the $N$-dimensional ball of radius $\rho$ and centre $0$
and
$C$ is a constant only depending on $N$, $D$, $q$, $\mu$, $\rho$, $T$, $p$ and
$\|u_0\|_{W^2_p(B_{\rho+2})}$ (and not on $n$). Set
$Q_\rho:=B_\rho\times(0,\rho)$.
Since $0\leq v_n\leq H$ too, the estimates yield
\[\begin{split}
\|v_n\|_{W^{2,1}_p(Q_\rho\times(0,T))}
&\leq C'
\left(1+\|v_n\|_{L^\infty(Q_{\rho+1}\times(0,T))}+
\|u_n\|_{W^{2,1}_p(B_{\rho+1}\times(0,T))}\right)\\
&\leq C'(1+H+CH),
\end{split}\]
with $C'$ only depending on $N$, $d$, $r$, $\mu$, $\rho$, $T$, $p$,
$\|f\|_\infty$,
$\|u_0\|_{W^2_p(B_{\rho+1})}$ and
$\|v_0\|_{W^2_p(Q_{\rho+1})}$.
This shows that the $(u_n)_n$ and $(v_n)_n$ are uniformly bounded in
$W^{2,1}_p(B_\rho\times(0,T))$ and $W^{2,1}_p(Q_\rho\times(0,T))$
respectively.

\noindent Step 4. {\em Existence of a solution.}\\
Now that we know that $(u_n)_n$ and $(v_n)_n$ are uniformly bounded in compact
sets with respect to the $W^{2,1}_p$ norm, taking $p>N+1$ and using the Morrey
inequality, we infer that this is also true with respect to the
$C^{\alpha}$ norm, for some $0<\alpha<1$. Then, by the Schauder estimates,
the time derivative and the space derivatives up to order $2$ are uniformly
H\"older continuous in compact sets too. As a
consequence, $((u_n,v_n))_n$ converges (up to subsequences) in
$C^{2,1}_{loc}$ to some $(u,v)$. Passing to the limit as $n\to\infty$ in
\eqref{un}, \eqref{vn} we eventually find that $(u,v)$ satisfies
\eqref{E}-\eqref{IC}. Form step 1 we know that $u$ and $v$ are bounded and
nonnegative. \hfill$\square$\\


\SE{The equation $h^L(c,\beta)=0$ }\label{sec:Rouche}

In this section we describe in detail how, for
$\xi>0$ small enough, equation \eqref{e4.8} admits two solutions close to
$\tau_+$ and $\tau_-$ respectively. We recall that $\tau_\pm=\pm
i\sqrt{(e/a)\xi}+O(\xi)$ are the roots of the trinomial
$g(\tau):=a\tau^2+d\xi\tau+e\xi$. Let us focus on $\tau_+$, the other case
being analogous.

Let $B$ be the ball of radius $A\xi$, centred at $\tau_+$; $A$ large and to be adjusted. 
For $\tau\in\partial B$, we have
$$|g(\tau)|=a|\tau-\tau_+|\, |\tau-\tau_-|
\geq aA\xi(|\tau_+-\tau_-|-
A\xi)=2aA\sqrt{e/a}\,\xi^{3/2}+O(\xi^2).$$
On the other hand we have
$|\varphi(\tau,\xi)|\leq C\xi^{3/2}+O(\xi^2)$.
We can therefore choose $A$ large enough and then $\xi$ small enough in such a
way that $|g|>|\varphi|$ on $\partial B$. Since $g$ and $\varphi$ are
holomorphic, by
Rouch\'e's theorem the equation \eqref{e4.8} has the same number of solutions in
$B$ as $g=0$, that is $1$. Notice that such a solution has positive imaginary
part proportional to $\sqrt{\xi}$ and real part of order $\xi$.

\end{appendix}

\section*{Acknowledgements}

This study was supported by the French "Agence Nationale de la Recherche"
through the project PREFERED (ANR 08-BLAN-0313). H.B. was also
supported by an NSF FRG grant DMS-1065979. L.R. was partially supported by the
Fondazione CaRiPaRo Project ``Nonlinear Partial
Differential Equations: models, analysis, and control-theoretic problems''.

\end{document}